\documentclass[12pt]{amsart}
\usepackage[english]{babel}
\usepackage{paralist,url,verbatim, anysize}
\usepackage{amscd}
\usepackage[all,cmtip]{xy} 
\usepackage{amsmath,amsthm,amssymb,amsfonts}
\usepackage[bookmarksopen=true]{hyperref}
\usepackage[usenames, dvipsnames]{color}
\usepackage{graphicx}
\usepackage{hyperref}
\usepackage{euler, amsfonts, amssymb, latexsym, epsfig,epic}

\makeatletter
\let\@fnsymbol\@arabic
\makeatother

\theoremstyle{plain}
\newtheorem{theorem}{\bf Theorem}[section]

\newtheorem{corollary}[theorem]{Corollary}
\newtheorem{lemma}[theorem]{Lemma}
\newtheorem{proposition}[theorem]{Proposition}

\theoremstyle{definition}

\newtheorem{example}[theorem]{Example}
\newtheorem{remark}[theorem]{Remark}

\newtheorem*{theorem*}{\bf Theorem}

\definecolor{mypink}{RGB}{215, 5, 234}

\begin{document}
\title{A universal tale of determinants and Gr\"obner bases}
\author{Aldo Conca}
\address{Dipartimento di Matematica, Dipartimento di Eccellenza 2023-2027, Universit\'a di Genova, Italy}
 \email{conca@dima.unige.it}
 
\thanks{The author was partially supported by  PRIN 2020355B8Y  ``Squarefree Gr\"obner degenerations, special varieties and related topics"  and by GNSAGA-INdAM. The author thanks the referees for  the careful reading of the paper and for the improvements of the exposition suggested.}
 

  \date{}
  
\dedicatory{Dedicated to Bernd Sturmfels on the occasion of his sixtieth birthday}

\begin{abstract}
In 1965 Buchberger defined  Gr\"obner bases  and an algorithm to compute them.   Despite a slow start,  already  in the eighties   Gr\"obner bases  had become  the main device  for symbolic computations involving polynomials as well as a theoretical tool for the investigation of ideals and varieties via the so-called Gr\"obner deformation techniques. 
Rings and algebraic varieties defined by means of determinants  are among the most classical objects in  commutative algebra, algebraic geometry  and invariant theory. By the end of the the eighties  the time was ripe for the computation of Gr\"obner bases of determinantal ideals. We will tell the tale of how (universal) Gr\"obner bases of determinantal ideals were identified and the key role played by Bernd Sturmfels and his collaborators in this enterprise. 
\end{abstract}
 
 \maketitle
 
 
\section{Generalities on Gr\"obner bases  and deformations}\label{Gbasis}
 
Let $K$ be a field and  $S$ be  the polynomial ring $K[x_1,\dots, x_n]$. As a $K$-vector space $S$ has very natural basis: the set of the monomials $\{ x^a : a\in {\mathbb N}^n\}$. 
A term order  on $S$  is a total order $\leq $ on the set of monomials of $S$ such that $x^a\leq x^b$ implies  $x^ax^c\leq x^bx^c$  and $1=x^0\leq x^c$ for all $c\in {\mathbb N}^n$. Given a term order $\leq $, any non-zero polynomial  $f\in S$ has an initial monomial  ${\operatorname{in}}(f)$, i.e.~the  largest monomial that appears in $f$. Furthermore any ideal $I$ of $S$ has an initial ideal ${\operatorname{in}}(I)=( {\operatorname{in}}(f) : f\in I\setminus\{0\} )$. A subset $G=\{f_1,\dots, f_s\}$ of non-zero elements of $I$ is a Gr\"obner basis of $I$ with respect to $\leq$ if ${\operatorname{in}}(f_1),\dots, {\operatorname{in}}(f_s)$ generates ${\operatorname{in}}(I)$ as an ideal. Equivalently $G$ is a Gr\"obner basis of $I$ if for every non-zero $f\in I$ there exists $f_i\in G$ such that ${\operatorname{in}}(f_i)$ divides ${\operatorname{in}}(f)$. The Noetherianity of $S$ implies that $I$ has a finite Gr\"obner basis and the Buchberger algorithm \cite[Chap.1]{BCRV} computes it  starting from any  finite set of generators of $I$. A universal Gr\"obner basis of $I$ is a set polynomials of $I$ that is a Gr\"obner basis with respect to all possible the term orders.  There are infinitely many term orders but, essentially,  only finitely many different Gr\"obner bases of $I$ so that $I$  has a finite  universal Gr\"obner basis. 
 
 The relationship  between $S/I$ to $S/\!{\operatorname{in}}(I)$ can be understood as the passage from the generic fibre to the special one of a $1$-dimensional flat family. An important consequence is that most of the  algebraic properties can be transferred form $S/\!{\operatorname{in}}(I)$ to $S/I$, see for example \cite{BC3,BCRV} for details.  For example $S/I$ is Cohen-Macaulay whenever  $S/{\operatorname{in}}(I)$ is Cohen-Macaulay.   By the very definition, ${\operatorname{in}}(I)$ is a monomial ideal and hence an object that can be investigated by means of combinatorial methods. The special case when ${\operatorname{in}}(I)$ is a radical ideal is very important for at least two reasons. The first is that a radical monomial ideal corresponds to a simplicial complex in the Stanley-Reisner duality and hence an entire toolkit of new ideas and resources can be employed in the investigation of the properties of $S/\!{\operatorname{in}}(I)$. The second is a recent result \cite{CV} asserting that when ${\operatorname{in}}(I)$ is radical the relationship between $S/\!{\operatorname{in}}(I)$ to $S/I$ is tighter, for example $S/\!{\operatorname{in}}(I)$ is Cohen-Macaulay if and only if  $S/I$ is Cohen-Macaulay. 
 
 \section{Generalities on determinantal ideals}\label{detring}
 Let $A$ be a matrix  with entries in a Noetherian ring $S$.  The ideal generated by all the minors (i.e.~ subdeterminants)  of $A$ of a given size, say  $t$,  is denoted by $I_t(A)$. It is called the determinantal ideal associated to $A$.  The classical Cauchy-Binet formula actually asserts  that $I_t(AB)\subseteq I_t(A)I_t(B)$ and this easily implies that $I_t(CAB)=I_t(A)$ whenever  $C,B$ are invertible matrices.  Said otherwise, if $f:F\to G$ is a $S$-linear map and $F$ and $G$ are finitely generated free $S$-modules the ideal $I_t(f)$ can be defined as the ideal of $t$-minors of any matrix representing $f$. 
 
 Determinantal ideals are ubiquitous in algebra and  geometry. Indeed many rings and algebraic varieties that are classically studied are, in one way or the other,  defined by determinantal conditions. 
 Among them a central role is played by the ``generic determinantal  ideals" that are defined as follows. 
 Let $K$ be a field  and $m,n$ be positive integers.  Let $S=K[X]=K[x_{ij} : (i,j)\in [m]\times [n]]$ be the polynomial ring in $mn$ variables and $X$  be $m\times n$ matrix whose $(i,j)$-th entry is equal to $x_{ij}$. 
 We call $X$ the generic $m\times n$ matrix (over $K$) and $I_t(X)$ the generic determinantal ideal. 
 The algebraic variety defined by $I_t(X)$ is the set of linear maps $f:K^n\to K^m$ of rank $<t$. For simplicity we will assume that $m\leq n$ which is not a restrictive assumption because ideals of minors do  not change under transposition. 

 The ideal $I_t(X)$ shows up  in classical invariant theory. Let us recall how. Consider the polynomial  ring $K[Y,Z]$ with $Y,Z$ generic matrices over $K$ of size $m\times (t-1)$ and $(t-1)\times n$. The group $G={\mathrm{GL}}_{t-1}(K)$ acts on  $K[Y,Z]$ by $Y\to YA$ and $Z\to A^{-1}Z$ for $A\in {\mathrm{GL}}_{t-1}(K)$.  Under mild assumption on $K$ (see \cite[Chap.7]{BV}) the ring of invariants $K[Y,Z]^G$ is generated by the entries of the product matrix $YZ$. Mapping $X\to YZ$ we get a surjective map $\phi:K[X]\to K[Y,Z]^G$. By rank reasons its  kernel clearly contains $I_t(X)$ and a result of classical invariant theory asserts that  $\ker \phi=I_t(X)$. Hence
 \begin{equation}
 \label{ringinv} 
 K[X]/I_t(X) \simeq  K[YZ]=K[Y,Z]^G.
 \end{equation} 
 In the early seventies Hochster conjectured that the ring of the invariants of any reasonably good group (i.e.~linearly reductive) acting on a polynomial ring over a field is Cohen-Macaulay. The conjecture was proved in full generality by Hochster and Roberts \cite{HR}  and reproved later on by Hochster and Huneke \cite{HH}  using  characteristic $p$ methods, in particular tight closure.  Other approaches to the Cohen-Macaulay property of $K[X]/I_t(X)$ have been developed more or less in the same years. One of them, which is important for our story, is based on  the so called standard monomial theory.  
 
   \section{Standard monomial theory}\label{stdmon}
  The polynomial ring $K[X]=K[x_{ij} : (i,j)\in [m]\times [n] ]$ has a very natural $K$-basis: the set ${\mathfrak{M}}$ of the monomials  in the $x_{ij}$'s. This basis is closed under multiplication and hence monomial ideals are (relatively) easy to deal with.  However to study the determinantal ideals $I_t(X)$ it would be better to have a $K$-basis ${\mathfrak{S}} $ of $K[X]$ such that $I_t(X)$  is generated, as a vector space, by a subset of the  ${\mathfrak{S}}$. Miraculously such a basis ${\mathfrak{S}}$ has been discovered in the seventies by Doubilet, Rota and Stein \cite[Thm.3.2.1]{BCRV} and its elements are certain products of minors of various sizes called standard monomials. The price to pay to deal with the standard monomial basis ${\mathfrak{S}}$ is that  it is not closed under product. But there is a replacement for this  drawback: the so-called straightening law that allows to keep under control the rewriting of the product of two elements in ${\mathfrak{S}}$  as a $K$-linear combination of elements of ${\mathfrak{S}}$. The entire package (the basis and rewriting rule) goes under the name of algebras with straightening law and it has been developed as an abstract theory in the seventies and eighties, see \cite{BV, E, DEP,MuS}. We will not enter into full details of this side of the story but we will only recall the definition of  Doubilet, Rota and Stein's standard monomials. A single $t$-minor $\delta$ of $X$ is denoted by $[r_1,\dots, r_t | c_1,\dots, c_t]$ where the $r_i$'s and the $c_i$'s are the row and column indices given in increasing order. 
The set of minors is partially ordered by 
$$[r_1,\dots, r_t | c_1,\dots, c_t]\leq [s_1,\dots, s_p | d_1,\dots, d_p] \mbox{  if and only if } 
\left\{ \begin{array}{l} 
 t\geq p, \\
 r_i\leq s_i \mbox{  for } i=1,\dots p\\
 c_i\leq d_i \mbox{ for } i=1,\dots p. 
\end{array} 
\right.
$$ 
 A standard monomial is a product of minors $\delta_1\delta_2\cdots \delta_v$ such that 
$\delta_1\leq \delta_2 \leq\dots \leq \delta_v.$ Writing the row  indices as rows of a Young tableau and similar for the column indices a standard monomial gets identified with a Young bitableau which is standard (i.e.~ with increasing rows from left to right and weakly increasing columns from top to bottom). For example, the standard monomial 
$$
[1, 3, 4, 5  |  1, 2, 3, 6]\cdot
[2, 3, 5 |  1, 4, 5]\cdot
[4 | 2] 
$$ 
is represented by the bitableau 
\begin{equation} 
\begin{array}{c}
\begin{tabular}{|c|c|c|c|c|c|} 
\hline
1 & 3 & 4  & 5 \\ \hline 
2 & 3 & 5 \\ \cline{1-3}
4   \\ \cline{1-1}
\end{tabular}\\
\end{array}
\begin{array}{c}
\begin{tabular}{|c|c|c|c|c|c|} 
\hline
1 & 2 & 3  & 6 \\ \hline 
1 & 4 & 5 \\ \cline{1-3}
2   \\ \cline{1-1}
\end{tabular}\\
\end{array}
 \end{equation}

\section{Gr\"obner bases of determinantal ideals via the RSK correspondence} 
 A term order on  $S=K[X]=K[x_{ij} :(i,j)\in [m]\times [n]  ]$ is called diagonal if the initial term of every minor $\delta$ is the product of the entries in the main diagonal, i.e.
 $${\operatorname{in}}( [r_1,r_2,\dots, r_t | c_1,c_2,\dots, c_t])=x_{r_1c_1}x_{r_2c_2}\cdots x_{r_tc_t}$$
 and antidiagonal if    
 $${\operatorname{in}}( [r_1,r_2,\dots, r_t | c_1,c_2,\dots, c_t])=x_{r_1c_t}x_{r_2c_{t-1}}\cdots x_{r_tc_1}.$$

 There are plenty of diagonal or antidiagonal term orders.  For example,  the lexicographic term order induced by 
$$x_{11}>x_{12}>\dots>x_{1n}>x_{21}>x_{22}>\dots>x_{mn}$$
is diagonal and  reverse lexicographic term order induced by the same total order of the variables is antidiagonal.    

\begin{theorem} 
\label{main} 
The set of $t$-minors of $X$ is a Gr\"obner basis of $I_t(X)$ with respect to every  diagonal or antidiagonal order. 
\end{theorem} 
Actually,  since what really matters are the initial terms of the minors and since diagonals and antidiagonals coincide up to an obvious symmetry, it is enough to prove the statement for a single diagonal or antidiagonal term order. 

This theorem has been proved, more or less  simultaneously,     by Sturmfels  \cite{S} and Caniglia, Guccione, Guccione \cite{CGG} at the end of the eighties. A version of  it  can be deduced from an earlier paper of Narashiman \cite{N}  which is unfortunately  written  using an unusual terminology and with no mention of Gr\"obner bases.   The proof given by Sturmfels  is a priceless jewel of beauty. We sketch it. 
 
The polynomial ring $K[X]=K[x_{ij} : (i,j)\in [m]\times [n] ]$ has hence two very natural $K$-bases: the set ${\mathfrak{M}}$ of the monomials  in the $x_{ij}$'s and the set of standard monomials ${\mathfrak{S}}$. Indeed they both consist of homogeneous elements.  Sturmfels observed that a version of the classical Robinson-Schensted-Knuth    correspondence (RSK for short, denoted by KRS in some papers) can be reinterpreted as an explicit, combinatorially defined and degree preserving bijection  between ${\mathfrak{S}}$ and ${\mathfrak{M}}$.  Hence it can be extended to a degree preserving $K$-linear map:
$${\mathrm{RSK}}: K[X]\to K[X]$$ 
sending each standard monomial to an ordinary monomial.  
Now since $I_t(X)$ has a $K$-basis of standard monomials its image  ${\mathrm{RSK}}(I_t(X))$ has the following two properties: 

\begin{itemize} 
\item[(a)] ${\mathrm{RSK}}(I_t(X))$  is generated by monomials,   
\item[(b)] ${\mathrm{RSK}}(I_t(X))$ and  $I_t(X)$ have the same Hilbert function. 
\end{itemize} 

 Note that (a) and (b) are typical properties of an initial ideal of $I_t(X)$.   However, by the very definition,  ${\mathrm{RSK}}(I_t(X))$ is only a $K$-vector subspace of $K[X]$ and not (yet) an ideal.  On the other hand, consider the ideal $J_t$ generated by the main diagonals of the $t$-minors  of $X$ and observe that if  
 \begin{equation}
 \label{schen1} 
 {\mathrm{RSK}}(I_t(X))\subseteq J_t
 \end{equation} 
 then we would have: 
  \begin{equation}
\label{schen2} 
 {\mathrm{RSK}}(I_t(X))\subseteq J_t \subseteq {\operatorname{in}}(I_t(X)).
 \end{equation} 
for every diagonal order.  Now in \eqref{schen2} the leftmost and rightmost vectors spaces are  graded and have the same Hilbert function.  Hence in \eqref{schen2} 
 one has equality throughout. In particular $J_t= {\operatorname{in}}(I_t(X))$, that is, Theorem \ref{main} holds and one has a bonus: 
 
 \begin{equation}
 \label{schen3} 
 {\mathrm{RSK}}(I_t(X))={\operatorname{in}}(I_t(X)).
 \end{equation} 
 
 The final observation of Sturmfels is that  the  inclusion \eqref{schen1} holds by a result of Schensted \cite{Sc} (again suitably reinterpreted) and this completes the proof  of Theorem \ref{main}.

 Another important  remark of Sturmfels is that,  for a diagonal term order,  the initial ideal ${\operatorname{in}}(I_t(X))$
 is the Stanley-Reisner ideal associated with a simplicial complex that can be identified with  the $(t-1)$-order complex of a  planar distributive lattice.  Bj\"orner  proved in \cite{Bj}  that such a simplicial complex is shellable and hence Cohen-Macaulay.  Summing, one establishes the Cohen-Macaulay property of the determinantal ring $K[X]/I_t(X)$ via Gr\"obner deformation.

   Sturmfels' approach to Gr\"obner bases of determinantal ideals and related ideas were  later extended  to other types of determinantal ideals and to their ordinary and symbolic powers, see   \cite{BC1, BC2, BC3,  C94, C98, HT}. In particular it turned out that  \eqref{schen3} holds for a very large class of determinantal ideals. 
   
 Back then it was  very surprising, and still is,  to discover that a combinatorial map, the RSK,  that has nothing to do with the selection of monomials in the support of polynomials actually ``computes" initial ideals and Gr\"obner bases of determinantal ideals.  It is another example of the unreasonable effectiveness of mathematical ideas.

   \section{ Universal Gr\"obner bases} 
   A universal Gr\"obner basis of an ideal $I$ in a polynomial ring is a subset of $I$ that is a Gr\"obner basis with respect to every term order. There are infinitely many term orders but,  as a single polynomial has  only finitely many initial terms, one can prove that every ideal $I$ has only a finite number initial ideals. This implies that $I$ has a finite universal Gr\"obner basis. 
 
 \subsection{Universal Gr\"obner bases for maximal minors}
 
 For the ideal of maximal minors $I_m(X)$ of the generic matrix $X$ the result is the best possible: 
 
 \begin{theorem} \label{BSZB} 
 The set of maximal minors of $X$ is a universal Gr\"obner basis of $I_m(X)$. 
\end{theorem} 

The original proof of Theorem \ref{BSZB}, obtained combining work of Sturmfels and Zelevinsky \cite{SZ} and of Bernstein and Zelevinsky \cite{BZ},  is based on a detailed analysis of the Newton polytope associated to the product of the maximal minors and on  the notion of matching fields. Indeed, in \cite{SZ} the proof of Theorem \ref{BSZB} is reduced to a  conjectural assertion on matching fields later  proved in \cite{BZ}. 
Another proof of  Theorem \ref{BSZB}  is given by Kalinin  in \cite{K}.  A simpler  and shorter proof of  Theorem \ref{BSZB} was later found of Conca, De Negri and Gorla \cite{CDG1}. It is based on the observation that under the substitution $x_{ij}\to a_{ij}x_{1j}$ 
 with $a_{ij}\in K$ generic elements  both $I_m(X)$ and the ideal generated by the initial terms of the maximal minors of $X$ equals the ideal of all the square-free monomials of degree $m$ in the variables 
 $x_{11}, x_{12},\dots, x_{1n}$.   This proof and the main results of \cite{C07, CS} were the starting points of a development that led to the identification and study of Cartwright-Sturmfels ideals, i.e.~multigraded ideals whose multigraded generic initial  ideal is radical, see \cite{CCC, CDG2,CDG3,CDG4, CDG5, CDG6, CW}.   A corollary of the proof  of Theorem \ref{BSZB} given in \cite{CDG1}  is the following: 
 
 \begin{corollary} 
 \label{UGBcorM}
 Let $D$ be any initial ideal of $I_m(X)$. Then: 
 \begin{enumerate}
 \item[(1)] $D$ is radical,
 \item[(2)] $S/D$ is Cohen-Macaulay,
 \item[(3)] $\beta_{ij}(S/D)=\beta_{ij}(S/I_m(X))$ for all $i,j$. 
 \end{enumerate} 
\end{corollary} 

 Corollary \ref{UGBcorM} (2) and (3) were  first proved by Boocher \cite{Bo} in a more general setting and follows also from the recent result of Conca and Varbaro \cite{CV} and the Cohen-Macaulayness of $K[X]/I_m(X)$.

\subsection{Universal Gr\"obner bases for $2$-minors}

For $t=2$ the isomorphism \eqref{ringinv}  takes the form: 

\begin{equation} 
\label{segre}
K[X]/I_2(X)\simeq K[ y_iz_j : (i,j)\in [m]\times [n]  ] \mbox{ where } x_{ij}\to y_iz_j. 
\end{equation} 
In other words,  $I_2(X)$ is the toric ideal, i.e.~a prime ideal generated by binomials,    associated with the classical  Segre variety.  Special techniques have been developed to study toric ideals and varieties, see the recent volumes  \cite{BG,CLS} for an overview.  In particular, Sturmfels wrote a fundamental  monograph on the Gr\"obner bases for toric ideals \cite{S2}. 
Given a finite set of vectors in  ${\mathcal{A}}=\{a_1,\dots, a_n\} \subset {\mathbb Z}^d$  one considers the associated toric ideal  $P_{\mathcal{A}}$, i.e.~the kernel of the map 
$$K[x_1,\dots, x_n] \to K[t_1^{\pm 1}, \dots t_d^{\pm 1} ] \quad x_i\to t^{a_i}.$$

In Chapter 4 of \cite{S2} Sturmfels defined: 

\begin{itemize} 
\item[(a)] the  set  $U_{\mathcal{A}}$ of  all the binomials that belong to some reduced Gr\"obner basis of $P_{\mathcal{A}}$, 
\item[(b)]  the set of circuits $C_{\mathcal{A}}$ where a circuit of $P_{\mathcal{A}}$ is  a binomial  in $P_{{\mathcal{A}}}$  that is  the  difference of coprime monomials  and has a  minimal support,  
\item[(c)] the Graver basis $G_{\mathcal{A}}$ of $P_{\mathcal{A}}$, i.e the set of the binomials $x^a-x^b$ in $P_{\mathcal{A}}$ such that there is no other binomial $x^c-x^d$  in $P_{{\mathcal{A}}}$ such that $x^c$ divides $x^a$ and $x^d$ divides $x^b$,
\end{itemize}
and observed that $U_{\mathcal{A}}$ is a universal Gr\"obner basis of $P_{\mathcal{A}}$. He also proved in  \cite[Prop. 4.11]{S2}: 

\begin{equation} 
\label{ugb2}
C_{\mathcal{A}}\subseteq  U_{\mathcal{A}}\subseteq G_{\mathcal{A}}
\end{equation} 

For special configurations ${\mathcal{A}}$, called unimodular, Sturmfels then proved in \cite[Prop. 8.11]{S2} that the equality 
\begin{equation} 
\label{ugb3}
C_{\mathcal{A}}=G_{\mathcal{A}}
\end{equation} 
holds. The conclusion is that $C_{\mathcal{A}}$ is a universal Gr\"obner basis of $P$ if ${\mathcal{A}}$ is unimodular. 

Returning to $I_2(X)$, taking the canonical bases $e_1,\dots,e_m$ and $f_1,\dots, f_n$ of ${\mathbb Z}^m$ and ${\mathbb Z}^n$ and 
$${\mathcal{A}}=\{ e_i+f_j :  (i,j)\in [m]\times [n] \} \subset {\mathbb Z}^m\oplus  {\mathbb Z}^n$$
the corresponding toric ideal $P_{\mathcal{A}}$ is exactly $I_2(X)$.  It turns out that  ${\mathcal{A}}$ is indeed unimodular. Hence one has a description of a universal Gr\"obner basis of $I_2(X)$ as the set of the circuits of ${\mathcal{A}}$.  The circuits,  in this case, are described as cycles in the complete bipartite graph $K_{m,n}$. In practice a circuit is described as follows. For a  pair
$(I, J )$ of sequences of distinct integers, say $I = i_1, \dots i_s$ , $J = j_1,\dots, j_s$  with $2\leq s \leq m$ and $i_k\in [m]$ and $j_k\in [n]$ for all $k$ one considers the binomial 
$$F_{I,J}=x_{i_1j_1}x_{i_2j_2}\cdots x_{i_sj_s}-x_{i_1j_2}x_{i_2j_3}\cdots x_{i_{s-1}j_{s}}x_{i_sj_1}.$$
It turns out that the set of the $F_{I,J}$ is indeed the set of the circuits of $I_2(X)$.  Summing up combining  \cite[Prop. 4.11]{S2} with \cite[Prop. 8.11]{S2}  one obtains  the following characterization of a universal Gr\"obner basis of $I_2(X)$.

\begin{theorem} \label{UGB2min} 
The set of   binomials $F_{I,J}$ is  a universal Gr\"obner basis of the ideal $I_2(X)$. 
\end{theorem}

 Theorem \ref{UGB2min} is  also a  special case of  a theorem proved by  Villarreal   \cite[10.1.11]{V} that describes  a universal Gr\"obner basis of the  toric ideal associated to the edge algebra of a graph.

\begin{remark} 
\label{minUGB2min}
Obviously  different pairs of sequences $I,J$'s  can give the same polynomial $F_{I,J}$ up to sign. 
A canonical choice can be made as follows: let $C$ be the set of polynomials $F_{I,J}$ where $I=i_1,\dots, i_s\in$ and $J=j_1,\dots, j_s$ are sequences of distinct elements such that: 
\[
i_1<i_k \mbox{ for } k=2,\dots, s\quad \mbox{and}\quad j_1<j_2.
\]
The set $C$ is indeed a  universal Gr\"obner basis of $I_2(X)$ and it is minimal, that is, no proper subset  is a universal Gr\"obner basis of $I_2(X)$. 
  \end{remark} 

\begin{example}\label{UGB3x3}
For a $3\times 3$ matrix $X$ the minimal universal Gr\"obner basis of $I_2(X)$ described in 
 Remark \ref{minUGB2min} consists of the $9$ minors of size $2$ and by $6$ binomials of degree $3$: 
\begin{gather*}
 F_{123,123}=x_{11}x_{22}x_{33}-x_{12}x_{23}x_{31},\quad 
 F_{123,132}=x_{11}x_{23}x_{32}-x_{13}x_{22}x_{31},\\
 F_{123,231}=x_{12}x_{23}x_{31}-x_{13}x_{21}x_{32},\quad 
 F_{132,123}=x_{11}x_{32}x_{23}-x_{12}x_{33}x_{21}, \\
 F_{132,132}=x_{11}x_{33}x_{22}-x_{13}x_{32}x_{21},\quad 
 F_{132,231}=x_{12}x_{33}x_{21}-x_{13}x_{31}x_{22}.
 \end{gather*} 
 \end{example} 
 
 An important consequence is the following. 

\begin{corollary} 
\label{UGBcor2}
Let $D$ be any initial ideal of $I_2(X)$. Then: 
\begin{enumerate}
\item[(1)] $D$ is radical,
\item[(2)]  $S/D$ is Cohen-Macaulay.
\end{enumerate} 
\end{corollary} 

Assertion  (1) follows immediately from Theorem \ref{UGB2min}. 
Statement (2)   can be seen either as a special case of the main result of \cite{CV}  or derived from the following general result of Sturmfels \cite[Ch. 8]{S2}.

\begin{theorem} \label{radinCM} 
Let  $<$ be a term  order and $I$ a toric ideal.  Then $\sqrt{ {\operatorname{in}}_<(I)}$ is the Stanley-Reisner ideal of a regular triangulation of a ball and hence it defines a Cohen-Macaulay ring.  
\end{theorem} 
 
\begin{remark}\label{SpecialIni}
The Betti numbers of the initial ideals of $I_2(X)$ are in general strictly larger than those of $I_2(X)$. For example if $m,n>2$  and $<$ a suitable term  order then the corresponding  initial ideal $D$ of $I_2(X)$  has a generator in degree $3$ (the initial  term of one of the binomials in Example \ref{UGB3x3}) so that $\beta_{13}(S/D)>\beta_{13}(S/I_2(X))=0$. 
Already for $m=n=4$ there is no initial ideal $D$ of $I_2(X)$ with $\beta_{2,j}(R/D)=\beta_{2,j}(R/I_2(X))$ for all $j$,  see Conca, Ho\c sten and Thomas \cite[Ex. 1.4]{CHT}. 
This behavior, (i.e.~no initial ideals of $I_2(X)$ have the Betti numbers of $I_2(X)$) can, at least for large $m$ and $n$, be  explained by a very interesting result of   Dao, Huneke and Schweig in \cite{DHS}. Indeed, if $I_2(X)$  has  an initial ideal $D$ with $\beta_{2,j}(S/D)=\beta_{2,j}(S/I_2(X))$ for all $j$ then $D$ is generated in degree $2$ and has only  linear first syzygies (because $I_2(X)$ has only linear first syzygies). For a monomial ideal $D$ generated in degree $2$ and with only linear syzygies  in \cite[Thm. 4.1]{DHS} it is proved that the regularity ${\operatorname{reg}} D$ of $D$ is bounded above by $\log_\alpha (f(mn))+\beta$ where the numbers $\alpha>1$ and $\beta$ and the linear function $f$ are given explicitly. It would then follow that 
$${\operatorname{reg}} I_2(X)\leq {\operatorname{reg}} D \leq \log_\alpha (f(mn))+\beta$$
 which is a contradiction since ${\operatorname{reg}} I_2(X)=m$ and $\log_\alpha (f(mn))+\beta$ growths much slower than that. 
\end{remark}

The $2$-minors do not form a universal Gr\"obner basis of $I_2(X)$ in general, nevertheless one has: 

\begin{proposition}\label{RLUGB2min}
The $2$-minors of $X$ are a universal revlex  Gr\"obner basis, i.e., they are a  Gr\"obner basis for all the reverse lexicographic orders. 
\end{proposition} 

This assertion, originally due to  Restuccia and Rinaldo \cite{RR}, can be proved by a careful analysis of the $S$-pairs reduction. It can be proved also with the help of a computer algebra system because, by the Buchberger algorithm, 
it is enough to prove it for a $3\times 3$ matrix  where that are ``only"  $9!$ reverse lexicographic orders and this number can be drastically reduced  by taking into consideration the obvious symmetries.

\subsection{Universal Gr\"obner bases for $t$-minors with $2<t<m$}  

A universal Gr\"obner basis of $I_t(X)$ is unknown for $2<t<m$.  Experiments with small values of $m,n,t$  show that $I_t(X)$ has initial ideals that are not radical and do not define Cohen-Macaulay rings.  
In particular:  
\begin{enumerate}
\item For $t=3$ and $m=n=4$ the ideal $I_t(X)$ has a  lexicographic initial ideal which is not radical. 
\item For $t=3$ and $m=4$, $n=5$ the ideal $I_t(X)$ has a lexicographic initial ideal $J$ such that  $J$ is not radical and  both $S/J$ and $S/\sqrt J$ are not Cohen-Macaulay. 
\item For $t=4$ and $m=n=5$ the ideal $I_t(X)$ has a  reverse lexicographic initial ideal which is not radical. 
\item For $t=4$ and $m=5$, $n=6$ the ideal $I_t(X)$ has a  reverse lexicographic initial ideal $J$ such that    $J$ is not radical  and both $S/J$ and $S/\sqrt J$ are not Cohen-Macaulay. 
\end{enumerate}

\begin{remark} 
For $t=3$  all the reverse lexicographic initial ideals we have computed are radical. For example, for  $m=n=4$ we have found $69$ distinct revlex initial ideals up to the natural  $S_4\times S_4$ and  ${\mathbb Z}/2{\mathbb Z}$ symmetry and they are all square-free.   It is hence possible that the ideal $I_3(X)$ has indeed a square-free universal revlex Gr\"obner basis.  This would indeed explain, from  a different perspective, a result proved by Conca and Welker  \cite{CW} asserting that, at least in characteristic $0$,   the ideal of $3$-minors of the  $m\times n$ matrix of variables remains radical after  setting  to $0$ an arbitrary subset of the variables.  
\end{remark}

   \section{Gr\"obner bases of determinantal ideals via secants} 
     Sturmfels and Sullivant published a paper \cite{SS} in 2006 proposing a new and simpler approach to the study of Gr\"obner bases of determinantal ideals. The  new insight  is to use  secant ideals and varieties.  Let us sketch briefly the main steps. The arguments work over any base field but, in order to avoid technical complications,  we assume right away that the base field is ${\mathbb C}$
 
 Let $I,J$ be  ideals  in a polynomial ring $S={\mathbb C}[x_1,\dots,x_n]$. 
The join  $I*J$ of $I$ and $J$ is by definition the kernel of the ${\mathbb C}$-algebra map 
\[S
\to S/I\otimes_{\mathbb C} S/J
\] 
that sends $x_i$ to  $y_i \otimes \ 1 + 1 \otimes  y_i $ for $i=1,\dots,n$.  Here $y_i$ denotes the coset of $x_i$ in the corresponding quotient ring.  More concretely, $I*J$ is obtained by eliminating the variables $y=y_1,\dots, y_n$ and $z=z_1,\dots, z_n$ from the ideal
$$I(y)+J(z)+( x_i-y_i-z_i : i=1,\dots, n)$$
where $I(y)$ is  the ideal obtained from $I$ by replacing the $x$'s with the $y$'s and similarly for $J(z)$. 
  
 If  $I$ and $J$ are homogeneous ideals defining projective varieties $X$ and $Y$ in ${\mathbb P}^{n-1}$ then the join  (ideal) $I*J$ defines the join variety $X*Y$ of $X$ and $Y$,  i.e.~the Zariski closure of the union of the lines   through  the points $p$ and $q$  with  $p\in X$ and $q\in Y$.   Furthermore, $I*J$ is radical if $I$ and $J$ are radical and prime  if $I$ and $J$ are prime.

More generally, given ideals $I_1,\dots, I_r$ of $S$ one defines the join ideal $I_1*I_2*\dots *I_r$ as the kernel of the map 
\[S
\to S/I_1\otimes_{\mathbb C} S/I_2\otimes_{\mathbb C}  \dots \otimes_{\mathbb C} S/I_r
\]
sending $x_i$ to 
\[
y_i \otimes 1 \otimes \dots  \otimes 1+1 \otimes y_i \otimes \dots  \otimes 1+\cdots +  1 \otimes 1 \otimes \dots  \otimes y_i.
\]

 The $r$-th secant ideal  of an ideal $I$ is the join $I*\cdots *I$ (with $r$ factors) and it will be denoted by  $I^{\{r\}}$.  Geometrically it defines the so-called $r$-secant variety, that is the closure of the union of the linear span of  $r$ points $p_1,\dots, p_r$ that belong to the variety defined by $I$.  The following two results,  proved in \cite{SS} and  \cite{SU},   are crucial for the application to determinantal ideals. 

\begin{lemma}\label{sec_in_monotone} 
Let $<$ be a  term  order on $S$ and $I$ be an ideal. 
Then  
\[
{\operatorname{in}}_<\bigl(I^{\{r\}}\bigr)\subseteq  {\operatorname{in}}_<(I)^{\{r\}}
\]
for all $r$.   
\end{lemma} 
 
  \begin{lemma}\label{SRsec} If $J$ is the Stanley-Reisner ideal of a simplicial complex $\Delta$ then 
  $J^{\{r\}}$ is the Stanley-Reisner ideal of the  simplicial complex $\Delta^{\{r\}}=\{ \cup_{i=1}^r  F_i : F_i\in \Delta\}$. 
  \end{lemma}

Back to the determinantal ideals $I_t(X)$ of the generic matrix $X=(x_{ij})$ and $S={\mathbb C}[x_{ij} ]$ equipped with an antidiagonal term order.  Let $J_t$ be the ideal generated by the  antidiagonals of the $t$-minors of $X$. 
By construction we have

\begin{equation}
\label{defi} 
J_t \subseteq  {\operatorname{in}}(I_t(X)). 
\end{equation}

Furthermore one has   
\begin{equation}
\label{seca1} 
I_{t}(X)= I_2(X)^{\{t-1\}}
\end{equation}
for $2\le t\le  m$. This equality can  be proved by purely algebraic arguments but  follows also   from the  well-known fact that any matrix of rank $\leq t-1$ is the sum of at most $t-1$ matrices of rank $\leq 1$.  

 Let $P=[m]\times [n]$ with poset structure defined by  $(i,j)\leq (h,k)$ if $i\leq h$ and $j\leq k$.  

By construction  $J_t$ is generated by the  antichains (pairwise incomparable elements) of $P$ of cardinality $t$.  In particular, $J_2$  is the Stanley-Reisner ideal of the order complex of $P$, i.e.~the simplicial complex $\Delta_2$ whose elements are the chains (totally ordered subsets) in  $P$. 
Denote by $\Delta_t$ the simplicial complex associated to $J_t$.  Since the union of $(t-1)$ chains cannot contain an antichain of cardinality $t$ and, because of  Lemma \ref{SRsec}, we have 
 $\Delta_2^{\{t-1\}}\subseteq \Delta_t$. On the other hand in this case it is easy to see that a subset of $P$  that does not contain an antichain of cardinality $t$ is actually  the union of  at most $t-1$ chains.  
 This is indeed what Dilworth's theorem states in general for posets:  the size of the largest antichain of any poset  equals the minimal number of chains needed to partition the poset. So we have $\Delta_2^{\{t-1\}}=\Delta_t$ or, equivalently, 
\begin{equation} 
\label{Qata}
J_t=J_2^{\{t-1\}}.
\end{equation} 
The final observation is that for an antidiagonal term order the $2$-minors of $X$ form a Gr\"obner basis,  i.e.  
\begin{equation} 
\label{2minok}
{\operatorname{in}}(I_2(X))=J_2.
\end{equation} 
 This can be checked easily directly using the Buchberger algorithm or using the toric presentation or  (for the lazy) using computer algebra system  as the ``action  takes place"  in a $3\times 3$ matrix at most.

Summing up if $<$ is an antidiagonal term oder on $S$, we have
$$J_t \underset{\eqref{defi}} \subseteq  {\operatorname{in}}(I_t(X))  \underset{\eqref{seca1}}  = {\operatorname{in}}(I_2(X)^{\{t-1\}}) \underset{\ref{sec_in_monotone}}  \subseteq  {\operatorname{in}}(I_2(X))^{\{t-1\}} \underset{\eqref{2minok}} =   J_2^{\{t-1\}} \underset{\eqref{Qata}}=  J_t.$$
  Therefore  one has equality throughout and  in particular  in \eqref{defi}. This proves Theorem \ref{main} for an antidiagonal order.  
  
 \section{Gr\"obner bases  via degrees and  multidegrees} 
 Another  approach to Theorem \ref{main} is presented in the book of Miller and Sturmfels \cite{MS} and it is based on ideas developed in the paper of Knutson and Miller \cite{KM}. To explain the main ideas, let us first recall the following   well known lemma.  Let $S=K[x_1,\dots, x_n]$ and $I$ be an homogeneous ideal. We denote by 
 $\dim S/I$ the Krull dimension of $S/I$,  by $\deg S/I$ its degree and  by $P_{S/I}(x)$ its  Hilbert polynomial.  With $e=\deg S/I$  and  $d=\dim S/I$  one has 
 $$P_{S/I}(x)=\frac{e}{(d-1)!} x^{d-1}+\mbox{ lower degree terms}.$$ 
Hence the knowledge of  $\dim S/I$ and $\deg S/I$ is equivalent to the  knowledge  of the largest degree component of $P_{S/I}(x)$.  
  
\begin{lemma}\label{kntrick}  Let  $I$ and $J$  be homogeneous ideals in $S$. 
 Assume that $J\subseteq I$, $\dim S/J=\dim
  S/I$, $\deg S/J = \deg S/I$ and $J$ is pure (i.e.~all its
  associated primes have the same dimension). Then $I=J$.
\end{lemma} 
 
Said otherwise, given an inclusion of homogeneous ideals  $J\subseteq I$ to prove equality one can try to prove the two ideals have the same Hilbert functions. But if $J$ is pure it is enough to prove that the highest degree component of the Hilbert polynomials of $S/J$ and $S/I$ agree.  Lemma \ref{kntrick} follows by simple considerations  on the primary decomposition and the additivity of the degree.

\begin{remark} 
\label{rema} 
Now suppose we have a homogeneous  ideal $I$  and we want to prove that  $f_1,\dots, f_s\in I$ are indeed a Gr\"obner basis of $I$ with respect to a term order. We look at $J=({\operatorname{in}}(f_1),\dots, {\operatorname{in}}(f_s))$. We have inclusion $J \subseteq {\operatorname{in}}(I)$ and we want to prove it is an equality. In view of Lemma \ref{kntrick} if $J$ is pure  we just need to prove that: 

\begin{itemize} 
\item[(1)]  $\dim S/J=\dim S/I,$
\item[(2)]  $\deg S/J = \deg S/I,$
\end{itemize} 
because then Lemma \ref{kntrick} applied to $J$ and ${\operatorname{in}}(I)$ implies that $J={\operatorname{in}}(I)$. Here we have used the fact that $\dim S/I=\dim S/\!{\operatorname{in}}(I)$ and  $\deg S/I=\deg S/\!{\operatorname{in}}(I)$. 
\end{remark} 

In the determinantal setting, with $S=K[x_{ij} : (i,j)\in [m]\times [n] ]$ we consider the ideal of minors  $I_t(X)$ and a diagonal or antidiagonal term order.  We apply  the scheme of reasoning of Remark \ref{rema}   to the ideal $I_t(X)$   to prove Theorem \ref{main}. We take $J_t$ to be the ideal generated by the  diagonals (or anti diagonals) of length $t$ in $X$.  The facets of the   simplicial complex $\Delta_t$ associated to $J_t$ can be describe as families of non-intersecting paths in the grid $[m]\times [n]$, see \cite{HT} for details. From this description it follows that   $J_t$  is pure and $\dim S/J_t=(m+n-t+1)(t-1)$.   Furthermore the degree of $S/J_t$ equals the number of facets of  $\Delta_t$. The later can be computed    using the Gessel-Viennot determinantal formula and one  obtains
 \begin{equation}
\label{HTdeg}
\deg S/J_t= \det \left( {m-i+n-j \choose m-i}
\right)_{i,j=1,\dots, t-1.}. 
\end{equation} 
 
Now the dimension of $S/I_t(X)$ is easily seen to be $(m+n-t+1)(t-1)$ and its degree have been computed by various authors (not using Gr\"obner bases), see  for example \cite{HaT, ACGH}. Therefore one can check   that degree and dimension of  $S/J_t$ and  $S/I_t(X)$ actually agree.  Some efforts are needed to check that the formulas given in \cite{HaT, ACGH} agree with \eqref{HTdeg} but usually one gets along by  using  Laplace expansions and Vandermonde determinants. Summing up  one obtains a proof of Theorem \ref{main} using Remark \ref{rema}.    

A critic to the argument just given is that it uses information about determinantal ideals (the formula for the degree) to compute a Gr\"obner basis, while, in many cases, e.g.~\cite{HT}, the goal is the opposite, i.e.~to compute Gr\"obner bases in order to derive information about the determinantal ideals as, for example their degree and Cohen-Macaulayness.

 The approach proposed by Knutson and Miller \cite{KM}  and Miller and Sturmfels \cite{MS} to Theorem \ref{main} does not use already known  results on the determinantal  rings and has the advantage to work  simultaneously  for a large class of determinantal ideals.  Another important point  is that a   larger class  is more suitable to inductive processes. The same type of strategy is used, under the name of  principal radical system, by  Hochster and Eagon \cite{HE}, to prove the Cohen-Macaulay property of determinantal rings. 
  
 In practice one identifies a carefully chosen class of  ideals, say $\{ I_w : w\in \Omega\}$ where $\Omega$ is a  set (or a poset),  such that $I_w$ is  generated by minors of $X$ of possibly different sizes  (that depend on $w$). Furthermore $I_t(X)$ is among the $I_w$'s.   For each $w\in \Omega $ one identifies  a somehow canonical set of generators of $I_w$ whose initial terms with respect to an appropriate term order generates a monomial ideal $J_w$.   To prove that canonical set of generators of $I_w$ is a Gr\"obner basis  using  Remark \ref{rema} one has to prove the following.

 \begin{lemma} For every  $w\in \Omega$ the ideal   $J_w$ is pure and  
 \label{lem1} 
 \begin{itemize}
 \item[(1)]  $\dim S/J_w=\dim S/I_w,$
\item[(2)]  $\deg S/J_w = \deg S/I_w.$
 \end{itemize} 
 \end{lemma}

As the minors of the matrix $X$ are naturally ${\mathbb Z}^m\times {\mathbb Z}^n$ multigraded both $S/I_w$ and $S/J_w$ acquire  a multigraded structure and hence have multigraded Hilbert polynomials and  multidegrees. 
Indeed, each finitely generated ${\mathbb Z}^m\times {\mathbb Z}^n$-module $M$ has an associated multidegree, denoted by  ${\operatorname{Deg}}(M)$,  which itself a homogeneous polynomial whose degree is the codimension of $M$, see \cite[8.5]{MS} for details. 
The   polynomial  ${\operatorname{Deg}}(M)$ captures both the codimension and, in finer form,  the degrees of $M$ in the various directions.  Summing up, to prove Lemma \ref{lem1} it is enough to prove:
 
 \begin{lemma}
 \label{lem2}  
For every $w\in \Omega$ the ideal $J_w$ is pure and 
 $${\operatorname{Deg}}(S/J_w)={\operatorname{Deg}}(S/I_w).$$
 \end{lemma} 
 
 The ideal $J_w$ is a generated by square-free monomials and the associated simplicial complex $\Delta_w$ has facets that can be described. In this way one checks that   $J_w$ is pure. Finally Lemma  \ref{lem2} is proved, roughly speaking,  by showing that both ${\operatorname{Deg}}(S/J_w)$ and ${\operatorname{Deg}}(S/I_w)$ satisfy the same recursion  involving  multidegrees of ideals of the same type for ``smaller" values of $ w$. This is usually a  hard task that requires  mastery in identifying and describing the ``right" combinatorial gadgets (e.g.~Schubert polynomials and pipe dreams in \cite{MS}). The final step is    to check the statement for the initial values of $w$ which is usually easy. 
 
 Yet another approach to Theorem \ref{main}, similar in spirit and inductive in nature,   makes use of the notion of linkage and has been applied by  Gorla, Nagel and Migliore to various classes of determinantal rings in \cite{GMN}. 
 
  \section{Conclusion} 
The search of Gr\"obner bases of determinantal ideals has been a very active research field in the last 30 years. Many important results have been achieved, many new techniques and ideas have been developed. In this process the contribution of  Bernd Sturmfels has been  fundamental.


\end{document}